\documentclass{article}[11pt]
\usepackage{epsfig}
\usepackage{epstopdf}
\usepackage{amssymb,amsmath}
\usepackage{graphicx}
\usepackage{caption}
\usepackage{apacite}
\usepackage{graphicx}
\usepackage{graphicx}
\usepackage{subcaption}
\usepackage[margin=1in]{geometry}
\usepackage[english]{babel}
\newtheorem{theorem}{Theorem}

\title{{\bf Labelling vertices to ensure adjacency coincides with disjointness}}

\author{Mahipal Jadeja, Rahul Muthu, Sunitha V.}
\begin{document}
\maketitle
\begin{quote}
{\bf Abstract} \small
Given a set of nonempty subsets of some universal set, their \textbf{intersection} graph is defined as the graph with one vertex for each set and two vertices are adjacent precisely when their representing sets have non-empty intersection. Sometimes these sets are finite, but in many well known examples like geometric graphs (including interval graphs) they are infinite. One can also study the reverse problem of expressing the vertices of a given graph as distinct sets in such a way that adjacency coincides with intersection of the corresponding sets. The sets are usually required to conform to some template, depending on the problem, to be either a finite set, or some geometric set like intervals, circles, discs, cubes etc. The problem of representing a graph as an intersection graph of sets was first introduced by Erdos $[1]$ and they looked at minimising the underlying universal set  necessary to represent any given graph. In that paper it was shown that the problem is NP complete.

In this paper we study a natural variant of this problem which is to consider graphs where vertices represent distinct sets and adjacency coincides with \textbf{disjointness}. Although this is really the same problem on the complement graph, for specific families of graphs this is a more natural way of viewing it. When taken across the spectrum of all graphs the two problems are evidently identical. Our motivation to look at disjointness instead of intersection is that several well known graphs like the Petersen graph and Knesser graphs are expressed in the latter method, and the complements of these families are not well studied. Thus our choice is justified and not merely an attempt to artificially deviate from existing work. We study the problem of expressing arbitrary graphs as disjointness graphs while minimising one of several possible parameters. The parameters we take into account are the maximum size of an individual vertex label, minimum universe size possible (disregarding individual label sizes) and versions of these two where it is required that each vertex gets labels of the same size. 
\end{quote}
\section{Introduction}\label{SecIntro}
Knesser graphs $KG_{n,k}$ are graphs whose vertices correspond to the
$k$ element subsets of an $n$ element set and two vertices are
adjacent precisely when their corresponding subsets are
disjoint. Clearly if $n<2k$ then the graph is an independent set of
vertices. If $n=2k$ then the graph is a matching. When $n=2k+1$ we get the
special family graph of {\bf odd graphs}. Knesser graphs are well studied. Many problems on them can be solved
clearly and efficiently using this set-theoretic definition. A natural
question, therefore, is to try and model an arbitrary graph in this
fashion. That is, come up with an underlying universal set and a choice
of unique subsets to associate with each vertex such that adjacency is
characterized by disjointness of the corresponding subsets. Clearly
for an arbitrary graph the above choice of all identical sized subsets
of a certain set will not work, because a graph defined in that manner
is necessarily vertex transitive.

This problem in arbitrary graphs can be viewed from several different
perspectives. It can either be a mathematical result
stating bounds or exact values of the size of these labels and the
universal set size or constructive methods providing algorithms for
them. In either case one may choose to try and minimise the universe
of labels used, independent of the individual label sizes and their
uniformity or non uniformity. One can also disregard the universe size
and minimise the size of each label in a uniform sized labelling. Towards this end we begin with a precise formulation of the problems as well as elementary results for well known classes of graphs.

The closely related concept is {\bf intersection graphs} for finite sets in which non-adjacency is characterized by disjointness of the corresponding subsets of underlying universal set. Here in our case, graph labels represent the intersection graph for  $\overline G$, not for G. So for a given graph, these two labelling approaches are entirely different (except for self-complementary graphs). Since the problems of intersection and disjointness on graph representation are equivalent, the disjointness version is also NP Complete. The problem of finding a vertex labelling for an arbitrary graph using distinct sets for different vertices  and all its standard variants we have listed are NP Complete. This follows from the fact that the equivalent problem of determining the intersection number of an arbitrary graph is NP Complete $[1], [7]$.

Intersection graphs have many applications in the
fields  of scheduling, biology, VLSI design and they are also used for development of faster algorithms for optimisation problems.

Firstly, Szpilrajn and Marczewski $[2]$ had observed that every graph is an intersection graph. For a graph with m edges and n vertices,  a  trivial upper bound for intersection number is m (see $[3]$). Tight upper bound of $n^2/4$ was proved by Erdos and he also proved a tighter bound for intersection number for graphs with strictly more than $n^2/4$ edges: p + t (see $[5]$), where p=number of pairs of vertices that are not connected by an edge and t is unique integer such that $t(t-1) \leq p < t(t + 1).$   Alon  Noga $[6]$ derived an upper bound as a function of maximum degree of a graph: $2e^2(d + 1)^2ln n$ where d=maximum degree of the complement graph of G and e=base of the natural logarithm.

In Section $2$ statements of problem variants are discussed along with some basic results. Results related to {\sc usn} are discussed in detail for some specific classes of graphs (including complement of complete graph, matching, paths, cycles, hypercube) in Section $3$. Effects of dynamic structural changes in graphs on {\sc usn} is discussed in Section $4$. Section $4$ also describes cartesian product based method to calculate {\sc usn}.  In Section $5$ and $6$, results are explained briefly for other variants, namely {\sc uusn} and {\sc iln} respectively. The final section concludes the paper. Future research directions related to our problem variants are also provided in this section.   


\section{Statements of Problem Variants}\label{probstat}
Here we define the variants of our problem which we consider.
{\bf Universe size number}({\sc USN}) of a graph is the number of elements in a
smallest universal set $S$ such that one can label the vertices of the
graph with distinct subsets of $S$ in such a way that adjacency
coincides with disjointness of the corresponding subsets. Empty set is not allowed to be used in labeling. 

{\bf Individual label number}({\sc ILN}) of a graph is the smallest size of the
largest label over all labellings of the vertices with uniques sets
such that adjacency coincides with disjointness. 

We also study variants of the above two problems where we require that all the labels have identical number of elements, referred to as the uniform variants. The parameters are thus denoted by {\sc uusn} and {\sc uiln} respectively.

\subsection{Results on {\sc usn}, {\sc uusn}, {\sc iln} and {\sc uiln}}
Here, we present results conecting the above parameters as well as general bounds on {\sc usn}.
\begin{theorem}
{\sc usn}$(G)$ $\ge$ $\lfloor{log_{2}n}\rfloor+1$ where $G$ has $n$ vertices. 
\end{theorem}
Using $\lfloor{\log_{2}{n}}\rfloor$ elements, at most $n-1$ non-empty subsets can be generated which can be assigned to at most $n-1$ labels of vertices. Here total number of vertices is $n$. In order to assign a non-empty as well as unique label to $n$th vertex $1$ additional element is required in underlying set of universe. Therefore value of {\sc usn} is at least $\lfloor{\log_{2}{n}}\rfloor+1$ for any given graph. 
\begin{theorem}
{\sc iln}$(G)\le$ {\sc usn}$(G)$.
\end{theorem}
Consider a valid labelling of any graph $G$, optimal  in terms of the universe size. Clearly, the number of elements used in total is {\sc usn}$(G)$. In this particular labelling, no vertex has more than {\sc usn}$(G)$ elements in its label, becasue the underlying set has {\sc usn}$(G)$ elements. This proves the result.
\begin{theorem}
{\sc uiln}$(G)=${\sc iln}$(G)$
\end{theorem}
Consider a valid labelling of any graph $G$, optimal in terms of {\sc iln}$(G)$. i.e. at least one vertex has {\sc iln} elements in its label. For each vertex $v$, add {\sc iln}$(G)-|label(v)|$ new elements to make all labels the same size. This process does not violate the requirements of disjointness based adjacency.
\begin{theorem}
{\sc uusn}$(G)\ge${\sc usn}$(G)$.
\end{theorem}
Compared to {\sc usn}, the {\sc uusn}  parameter has an additional constraint of uniformity of individual vertex label sizes. Therefore the value of {\sc uusn} is at least {\sc usn}.

\begin{theorem}
{\sc usn}$(G+v)=${\sc usn}$(G)+1$, if $d(v)=n(G)$.
\end{theorem}
Since the vertex $v$ is adjacent to all other vertices, no element used in its label can be used in the labels of any other vertex. Disregarding $v$, the graph $G$ requires {\sc usn}$(G)$ elements for labelling its vertices even as a subgraph of $G'$. There is no purpose served in having more than one element in the label of $v$ because $v$ has no non-neighbours.

As a corollary we have the following theorem.
\begin{theorem}
{\sc usn}$(K_n)=n$
\end{theorem}

%
%
%
%
%
  
\section{Results on {\sc usn} for some special families of graphs}

-In this section, we dreive results on {\sc usn} (either exact or asymptotic) on the following classes of graphs: Paths, Cycles, Wheel Graph, Hypercube, Complete Graph, Complement of Complete Graph $(\overline{K_n})$, Matching $(M_{2n})$).


\begin{theorem}
For matching, ${\text{\sc usn}}=O(\log_{2}{2n})$

\end{theorem}
Consider matching graph on $2n$ vertices. (Number of edges=$n$). 
In a matching each vertex of the graph is adjacent to exactly one other vertex. Hence, the label of each vertex must have non-empty intersection with the labels of all the remaining $2n-2$ vertices.
Assume that the underlying universal set $U$ has $k$ elements. Our aim is to calculate the value of $k$.
We can summarize our requirements:
\begin{enumerate}
	\item $|k|$ should be large enough to generate at least $2n$ of non-empty subsets. 
	\item  Each subset must have non-empty intersection with all remaining $2n-2$ subsets.
\end{enumerate}
If we consider only subsets of cardinality $k/2$ then each subset will have exactly one disjoint subset and it will have non-empty intersection with all the remaining subsets. (Pigeonhole principle)
In order to fulfill the first requirement, ${{{k} \choose {k/2}} \geq 2n}$.
By Stirling's approximation, {\sc usn}=$k$=$O(\log_{2}{2n})$
\begin{theorem}
{\sc usn}$(\overline{K_n})=1+\lceil\log_2{n}\rceil$
\end{theorem}
Proof:\\
Case I:\\ 
Consider  $n=2^k$,\newline $k\in \mathcal{N}$ and $A=\{1,2,3,\ldots,k\}$\\
Clearly, $A$ has $2^k$ distinct subsets, out of which only $2^k -1$ subsets can be used for labeling, since empty set is not allowed. \newline
Using all non-empty subsets of $A$ assign labels to all $n-1$ vertices of the graph. 
Assign $\{k+1\}$ as a label to $n^{th}$ vertex (because all non empty subsets of $A$ are already used).
Add $\{k+1\}$ to all labels of the graph excluding the $n^{th}$ vertex. 
All vertices are pairwise disjoint and their labels contain the common element $\{k+1\}$, so intersection of any non-adjacent vertex pairs is non-empty. Hence the final labelling is valid. 
Hence an optimal labeling using $k+1$ is possible, i.e.$1+\log_2{n}$ elements.\\
Therefore {\sc usn}$=1+\log_2{n}$
\newline Case $II$:\newline
{\sc usn} is at least $k+2$.where number of vertices:$n=2^k+x$ and $x \geq 1$ and $x < 2^k$\newline \newline
Proof for lower bound:\\
Consider $n=2^k+x$,\newline $k \in \mathcal{N}$ and $A=\{1,2,3,\ldots,k,k+1\}$\\
Observation:\newline 
Here disjoint subsets are not allowed for labeling of independent set.\\
So at most $1$ subset can be used for labeling from each pair $(S,A-S)$.\\So the total number of available sets for labeling is reduced by half.\\
Here $|A|=\{k+1\}$. So total number of subsets are $2^{k+1}$.\\ Out of these $2^{k+1}$ subsets only half i.e. $2^{k}$ (at most) can be used for labeling. But here the number of vertices: $n=2^{k}+x$ and $x\geq 1$. Number of vertices are at least $2^k+1$.\\Hence, it is impossible to label $2^k+1$
vertices with $2^k$ subsets.\\ So at least one extra element is required in set $A$.\\ 
Hence the lower bound is $k+2$.\\
Assign valid and optimal labels to any  subset of $n - x$ vertices using the previous method. Select any $x$ labels and take the complement of these labels with respect to the set $A \cup \{k+2\}.$
Add $\{k+1\}$ to all the $x$ newly created labels.  
Assign these $x$ new labels to the remaining $x$ vertices.
Note: Here all labels are unique i.e. no repetition of label and common element present in all labels is $\{k+1\}$.\\
So the final labeling is valid and optimal.
\begin{theorem}
{\text{\sc usn}}  $K_{s,t}= 2+\lceil\log_2s\rceil +\lceil\log_2t\rceil$\
\end{theorem}
	Complete bipartite graph consists of two independent sets.\\From the valid and optimal labeling of one independent set nothing can be used in the second independent set. So labeling of these two independent sets should be entirely disjoint. So exactly $1+\lceil\log_2s\rceil +1+\lceil\log_2t\rceil$=$2+\lceil\log_2s\rceil +\lceil\log_2t\rceil$ elements are required in the underlying universe set.
\begin{theorem}
${\text{\sc usn}}(P_n) = O(logn)$.
\end{theorem}
\begin{figure}[htb]
\centering
  \begin{tabular}{@{}cccc@{}}
    \includegraphics[width=.450\textwidth]{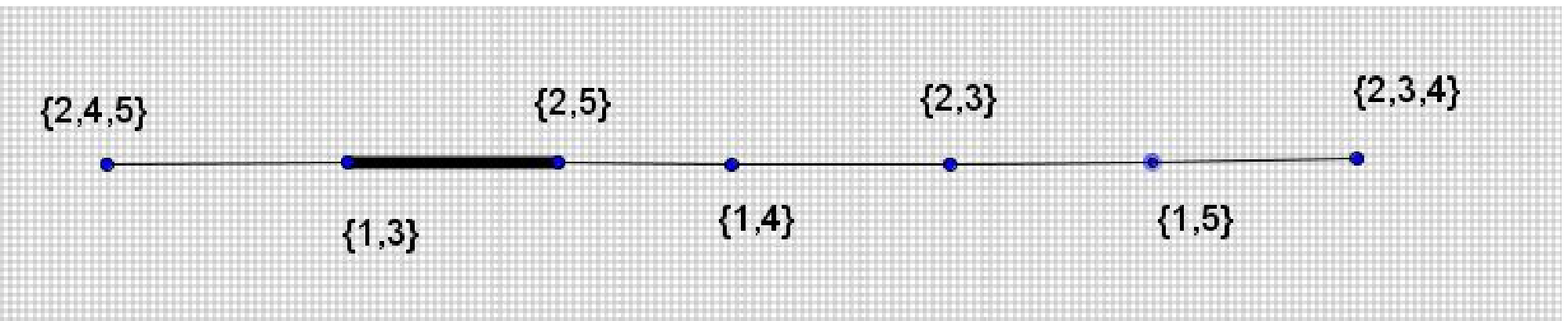} &
    \includegraphics[width=.450\textwidth]{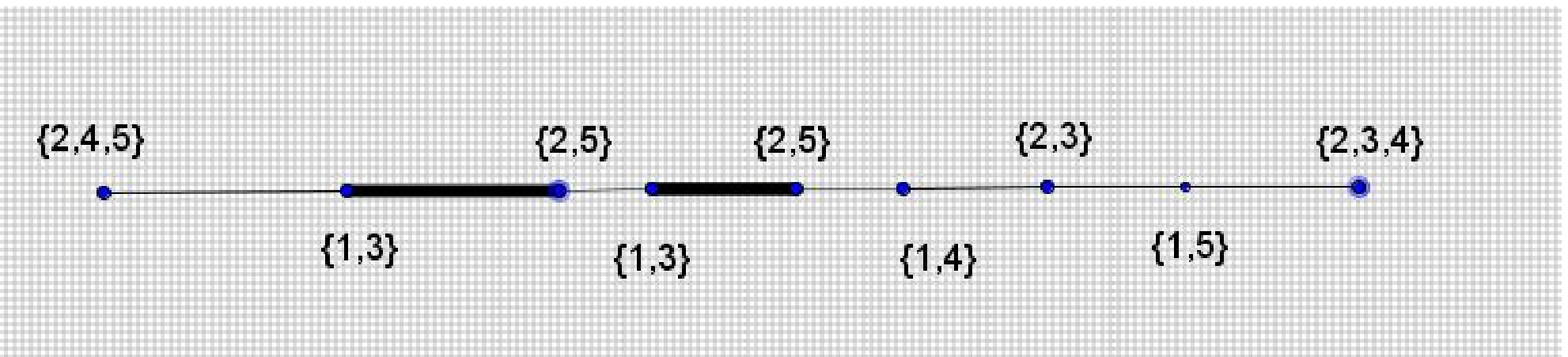} \\
		\includegraphics[width=.450\textwidth]{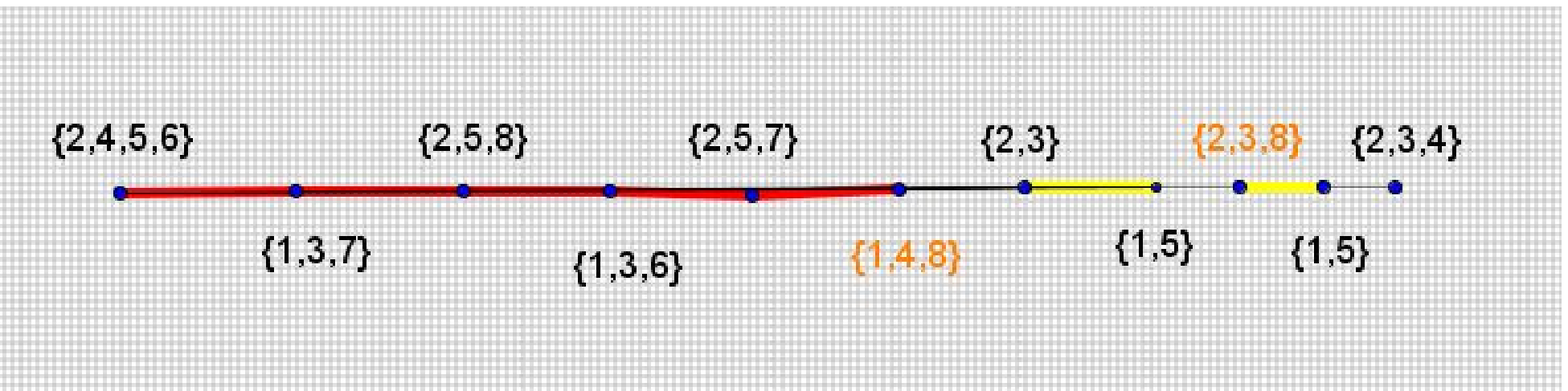} &

    \multicolumn{2}{c}{\includegraphics[width=.45\textwidth]{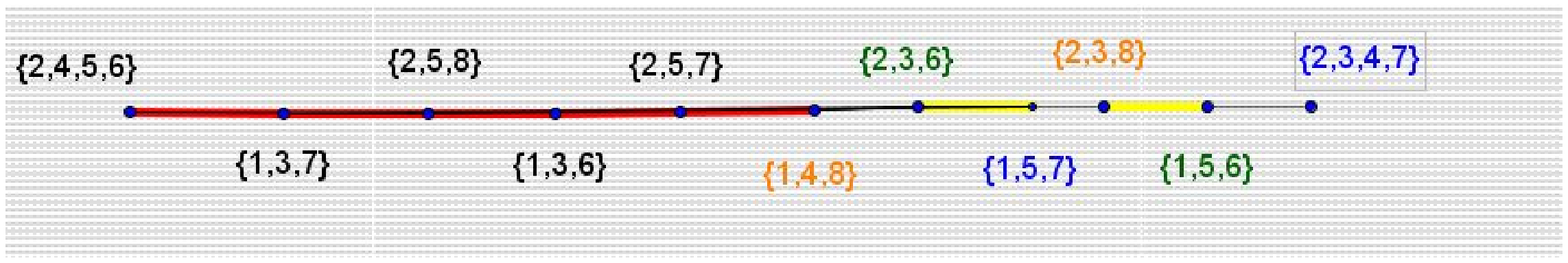}}
  \end{tabular}
  \caption{ Addedge procedure : $P_{7}$ to $P_{11}$ }
\end{figure}

Input: Consider an optimal and/or valid labelling of any given path ($P_{n}$) with universe size number $k$.\newline
Output:Valid labelling of  $P_{n+2}$, $P_{n+4}$, $P_{n+6}\ldots...,P_{5n/3}$ with {\sc usn}=$k+3$.\\
\newline
Procedure Addedge: \\
Before adding an edge :\\

Step 1:\\
Identify 4 consecutive vertices. 
Say valid labels for these vertices $v_{1}$, $v_{2}$, $v_{3}$ and $v_{4}$ are $m, p, q$ and $r$ respectively.\\
After adding an edge with 2 new vertices $v_{2'}$ and $v_{3'}$:\newline 
Step 2:\\ 
 p' and q' are labels of new vertices $v_{2'}$ and $v_{3'}$ respectively.\\
Assign $p'=p$ and $q'=q$.\newline
Observation After step 2 :\\
All vertices except $v_{1}$, $v_{2}$, $v_{3}$ and $v_{4}$ will preserve non-adjacency with  $v_{2'}$ and $v_{3'}$ because $v_{2}$, $v_{3}$ and $v_{2'}$, $v_{3'}$ have common non-neighbors.
\newline Step 3(a):\newline  
$v_{1}$ is not adjacent to $v_{2'}$. \\
$p'=p' \cup \{a\}$
$m=m \cup \{a\}$ where $a \notin {\text {\sc usn}}$ \newline

Step 3(b):\newline
$v_{4}$ is now non adjacent to $v_{3}$. \\
$q=q \cup \{b\}$\\
$r=r\cup \{b\}$ where $b \notin {\text{\sc usn}}$ and $b\neq a$
\newline step 3(c) :\newline
$p$ and $p'$ are non adjacent and they have common elements because $p\subseteq p'$.
$q$ and $q'$ are non adjacent and they have common elements because $q'\subseteq q$.
But $p$ and $q'$ are non-adjacent and they don’t have common element. So add one element $c$ to both of them.\\ 
i.e. $p=p\cup \{c\}$\\
     $q'=q'\cup \{c\}$ where $c \notin {\text{\sc usn}}$, $c\neq b$ and $c\neq a$.\\
	  So after step 3c) cardinality of $m, p, q, p',q'$ and $r'$ is increased by exactly 1 using 3 new elements namely $a, b$ and $c$. \\
New {\sc usn}=old {\sc usn}+3.
\newline
Upper bound calculation:\\
For any given  path graph on n vertices with {\sc usn} k, 
Addedge procedure can be applied at most $n/3$ times to get {\sc usn} of size k+3.
(edge ($v_{1}$, $v_{2}$) and ($v_{3}$, $v_{4}$)  cant be used for Addedge because use of any of them will increase {\sc usn} by 1 more.)\\
So by using one edge Addedge will eliminate possible use of 2 other edges.\\ So at most $\left\lfloor\frac{n}{3}\right\rfloor$  edges can be used and they will generate at most $2n/3$ new vertices.\\
Recurrance: $T(5n/3)=T(n)+ 3$\\            $T(n)=T(3n/5)+3 $\\
{\sc usn} : O$(\log{n})$ \\
Note: We can generate paths for all values of $n$ by applying the addedge procedure repeatedly. One can generate all $P_{2x+1}$  by repeatedly applying the addedge procedure on a valid and optimal labelling of $P_{7}$, similarly all $P_{2x}$ can be generated by application of the addedge procedure on a valid and optimal labelling of $P_{6}$. Here $x > 3$ and $x \in \mathcal{N}$.
\begin{theorem}
${\text{\sc usn}}(C_n)= O(\log{n})$.
\end{theorem}
For cycles results are slightly better for some values of $n$ because $C_{n}$ contains exactly one more edge as compared to $P_{n}$ and hence it is possible to apply addedge procedure exactly $1$ more time as compared to $P_{n}$ (specifically  when n is multiple of $4$).
\begin{theorem}
${\text{\sc usn}}(W_n) = O(\log{n})$.
\end{theorem}
Wheel graph consists of $C_{n-1}$ and one additional vertex with degree $n-1$. Therefore, for all parameters parameter value for $W_{n}= \text{value for }  C_{n-1}+1$. Because reusing of any labels is not possible for the vertex with degree $n-1$. In order to assign non-empty unique label, one extra element is required. 
\begin{theorem}
${\text{\sc usn}}(Q_n) < 3n+O(logn)$
\end{theorem}
Proof:
For hypercube, the valid labelling problem can be viewed as two independent subproblems:

\begin{enumerate}
	\item  Valid labelling of layers in order to preserve non-adjacency between them. (except adjacent layers-if at least one edge is present between two layers then those two layers are adjacent).
i.e. non-adjacent layers should have at least one element in common to preserve non-adjacency.

\item Valid labelling to preserve adjacency and non-adjacency between adjacent layers. 
\end{enumerate}
{\bf Algorithm:}\\
Input: $Q_{n}$\\
Output: Valid labelling of $Q_{n}$ with at most $3n+O(\log{n})$ labels.\\
\newline 
Step 1:\\
Consider $P_{n+1}$( Path of length $n$) for given $Q_{n}$\\
Apply path-labelling algorithm to get a valid labelling of $P_{n+1}$ with 
O$(\log{n})$ labels.\\
Now, For all $i$, Labels of all vertices at layer$i$= Label of vertex $P_{i+1}$
For example, for $Q_{6}$ consider a valid and optimal labelling of $P_{7}$ and
 assign 245, 13, 25, 14, 23, 15, 234 to the 7 layers of $Q_{6}$. 
So far, O$(\log{n})$ total number of labels are used and after this step all non-adjacent layers will preserve non-adjacency.\\
Additionally all vertices within one layer will have at least one common element. So non-adjacency within one layer will be preserved too.\newline
Step 2: \\

(a)Consider three disjoint (disjoint with used labels of step 1 too) sets each of size n. $S_{1}$, $S_{2}$ and $S_{3}$.\\
Now set $S_{1}$ will be used for labelling of layers: 3, 6, 9...\\
Similarly set $S_{2}$ will be used for labelling of layers: 1, 4, 7, 10...\\
and set $S_{3}$  will be used for labelling of layers: 2, 5, 8, 11....\\
Now $Q_{n}$ contains $n \choose i$ vertices at layer i.\\ 
Assign $n \choose i$ different labels to all vertices of layer $i$ using the $n \choose i$ subsets of  the underlying set $S$ mentioned above. (Using conventional subset assignment method to hypercube i.e. layer $i$ should contain all subsets of length $i$ ) \\
So  after this step 2 (a), all labels will be unique.\\
Note:For $0^{th}$ layer use valid labelling of $P_{1}$ vertex obtained from step 1. \\
After step 2(a), adjacency is preserved fully because all adjacent vertices have disjoint labels because they are using disjoint underlying sets for labelling.\\
Total $3n+O(\log{n})$ labels are used so far.\\
(b): \\
Now non-adjacency between two adjacent layers is remaining to be considered. 
For each vertex of layer $i$, add extra elements  from layer $i-1$ to preserve non-adjacency with layer $i-1$. \\
Procedure to add extra elements to label  of vertex $v$ of layer $i$ :\\
\begin{enumerate}
	\item Consider all vertices which are adjacent to vertex $v$ and at layer $i-1$. Take union of all labels of these vertices, say $P$ \\
  \item Take the complement of this set $R$ with respect to the underlying labelling set of layer $i$. (one of $S_{1}$, $S_{2}$ and $S_{3}$) say $Q= P^{c}$ \\
\item  New label of $v$= old label of $v \cup {Q}$
\end{enumerate}
Now the $i^{th}$ layer will contain elements from both sets: set used for layer $i$ and $i-1$.\\

Add extra elements only from the originally assigned underlying set in 2(a). i.e. for layer $i+1$ add extra elements to preserve non-adjacency from the underlying set of layer $i$ not $i-1$.\\
Final labelling after step 2(b) will preserve adjacency as well as non-adjacency for all vertices. 

\section{Effect of dynamic graph operations on {\sc usn}}
In this section we study the effect of addition/deletion of a vertex/edge on the {\sc usn} of a graph. We also derive an upper bound (though not tight) on the {\sc usn} of the union of any two graphs on the same vertex set, in terms of their individual {\sc usn} values.

\begin{theorem}
A vertex can use a singleton label  in some valid labelling, if and only if the subgraph induced by its non-neighbours has zero edges. 
\end{theorem}
Suppose a vertex $v$ has two non neighbours $u$ and $w$ such that $(u,w)\in E(G)$. Suppose the labelling of $v$ is a singleton set. Now due to nonadjacency to $v$, both $u$ and $w$ must contain that element, but due to their adjacency to each other, $u$ and $w$ cannot contain a common element. The converse is trivial. This proves the lemma.
\begin{theorem}
{\sc usn}$(G+v)=${\sc usn}$(G)+1$, if $d(v)=n(G)$.
\end{theorem}
We have stated and proved this result in Section $2$ and have merely restated it here for the sake of completeness. 

\begin{theorem}
{ $G'=G+{v}$ where order of G is n.\\ ${\text{\sc usn}}$(G)$ \le {\text{\sc usn}}(G') \le {\text{\sc usn}}$(G)$+(n-1)$}
\end{theorem}
\begin{theorem}
{ $G'=G-{v}$ where order of G is n.\\ ${\text{\sc USN}}$(G)$-(n-1) \le {\text{\sc USN}}(G') \le ${\text{\sc USN}}$(G)$}
\end{theorem}
We give, here, a combined proof of the previous two results. Adding(removal) a vertex can not decrease(increase) the value of {\sc usn}.\\In the worst case the existing graph $G$ may have clique of size $n-1$ and one isolated vertex, say $v_i$ with $|label|=${\sc usn}. Consider the case where the newly added vertex $v$ is adjacent to only $v_i$. $v$ is non-adjacent to the clique of size $n-1$ and in order to preserve non-adjacency with the clique extra $n-1$ elements are required.   
\begin{theorem}
$G'=G\setminus e$, where $e\in E(G)$. 
{\sc usn}$(G')\le${\sc usn}$(G)+1$.
\end{theorem}



\begin{theorem}
$G'=G+{e}$, where $e\notin E(G)$. 
{\sc usn}$(G')\ge{\text{\sc usn}}(G)$-1.
\end{theorem}
After Removal of an edge, {\sc usn} may increase  by at most one. This is because in the worst case, only one new element need be added in the labels of the endpoints in order to respect non-adjacency.\\ This proves Theorem 18, and Theorem 19 is effectively just a rewording of Theorem 18.








We now present an algorithm to compute a valid labelling of any graph. 

{\bf Algorithm to find a valid labeling for any given graph\\}
{\bf Step 1:}For any  graph  on $n$   vertices   start with the optimal valid labelling of the  corresponding  complete  graph  $K_n$  {\sc usn} $=n$.\\
{\bf Step 2:} Delete the necessary set of edges one by one from this $K_n$ to transform complete graph into the given graph. After deleting each edge, add an extra element to labels of both endpoints of that particular edge.\\
This algorithm will give valid labeling of any given graph with at most $\sum n$ elements.\\   

The next set of results gives an upper bound on the {\sc usn} of the union of two graphs on the same vertex set, in terms of the individual {\sc usn} values.
\begin{theorem}
Let $A, B, C, D$ be sets. Then $(A \times B) \cap (C \times D) \neq \emptyset$ if and only if $A\cap C \neq \emptyset$ and $B\cap D \neq \emptyset$.
\end{theorem}
\begin{theorem}
Let $G$ and $H$ be two graphs on the same vertex set $V$. Further, suppose $E(G)\cap E(H)=\emptyset$. Then {\sc usn}$(G+H)\le$ {\sc usn}$(G) \times ${\sc usn}$(H)$.
\end{theorem}
Take optimal labellings of the vertex set using disjoint universes for the graphs $G$ and $H$. Now consider the graph $G+H$. The vertex sets of $G$ and $H$ are identical. For each vertex $v$ in $G+H$ give it the label $l_G(v) \times  l_H(v)$. Clearly two vertices are nonadjacent only if they are nonadjacent in both $G$ and $H$.  In that case their labels under the two labellings will each be intersecting. From Theorem 20, it follows that their cartesian product new label will also intersect. Similarly for the case of non-intersection (adjacent vertices).
\begin{theorem}
From the above result we can infer that the {\sc usn}$(P_n)\le (1+\log \frac{n}{2})^2$.
\end{theorem}
The path is the union of two disjoint matchings. Each matching has {\sc usn} $O(\log n)$. From Theorem 21, we see that the graph has {\sc usn} $O (\log n)^2$. We state this here just as an application since we have a better bound for paths (Theorem 10).
\section{Results on {\sc uusn} }
\begin{theorem}
{\sc uusn}$(P_n) = O(logn)$.
\end{theorem}
Modified addedge procedure\\
Input: Consider an optimal and/or valid uniform labelling of any given path ($P_{n}$) with uniform universe size number $k$.\newline
Output:A valid uniform labelling of  $P_{n+2}$, $P_{n+4}$, $P_{n+6}\ldots P_{5n/3}$ with {\sc uusn}=$k+3$.\\

After $i^{th}$ iteration of the addedge procedure, the cardinality of the labels of the first $3i+1$ vertices will be increased by exactly 1. So the first $3i+1$ vertices of the path has a uniform labelling. The label of the $(3i+1)^{th}$ vertex will certainly contain exactly one of the 3 elements namely $a,b,c$. WLOG $a$ is used in the label of the $3i+1^{th}$ vertex.\\ 
Partition the remaining vertices $V\setminus\{V_{1},V_{2},\ldots,V_{3i+1}\}$ into $2$ sets: $A$: Even numbered vertices and $B$: Odd numbered vertices.\\
Add $b$ to all labels of $A$ and add $c$ to labels of $B$. After this step, the cardinality of the remaining vertices are also increased by $1$.   
In general, the cardinality of each vertex is increased by exactly $1$. 
So the final labelling is uniform and valid. 
No additional elements are required in this procedure.
Hence {\sc uusn}$=O(\log{n})$
By applying the modified addedge procedure, the following two results are derived. 
\begin{theorem}
{\sc uusn}$(C_n) = O(logn)$.
\end{theorem}
\begin{theorem}
{\sc uusn}$(W_n) = O(logn)$.
\end{theorem}
\begin{theorem}
{\sc uusn}$(Q_n) < 3n+O(logn)$
\end{theorem}
The hypercube labelling algorithm, adds exactly $n$ elements in each label (except for labels of the first two layers of the hypercube). 
If the underlying path on $n+1$ vertices has a valid uniform labelling with $O(\log{n})$ elements then it is possible to obtain a final uniform labelling using $3n+O(\log{n})$ labels using following modified algorithm.
\begin{itemize}
	\item Apply hypercube algorithm to get a valid non-uniform labelling (non uniform just because of the first two layers).
	\item The first layer contains only a single vertex say $v$. New label of $v= \text{old label  of v }  \cup S_{k}$, where $S_{k}$ is underlying set for labelling which is used in the $3^{rd}$ layer of the hypercube.
	\item Labels of the second layer contains only one additional element apart from the elements of the $2^{nd}$ vertex of the corresponding path. Consider $(n-1)^{th}$ layer of the hypercube which contains all $n-1$ sized subsets of the underlying labelling set. Layer 2 and $n-1$ both contains the same number of elements. Add all $n-1$ element subsets generated at layer $n-1$ into their corresponding copy in the $2^{nd}$ layer.
\end{itemize}
\begin{figure}[htb]
\centering
  \begin{tabular}{@{}cccc@{}}
    \includegraphics[width=.50\textwidth]{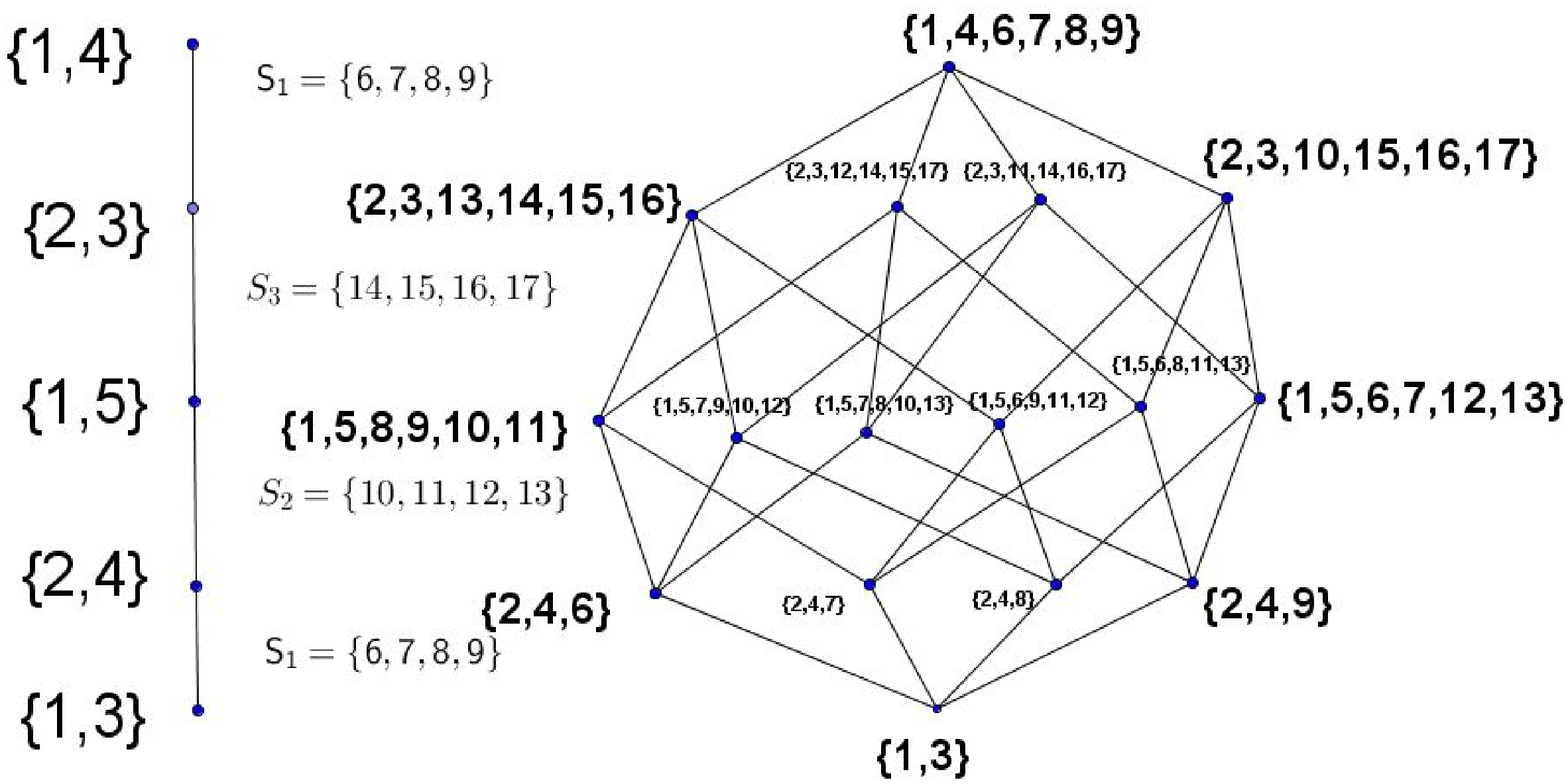} &
    \multicolumn{2}{c}{\includegraphics[width=.50\textwidth]{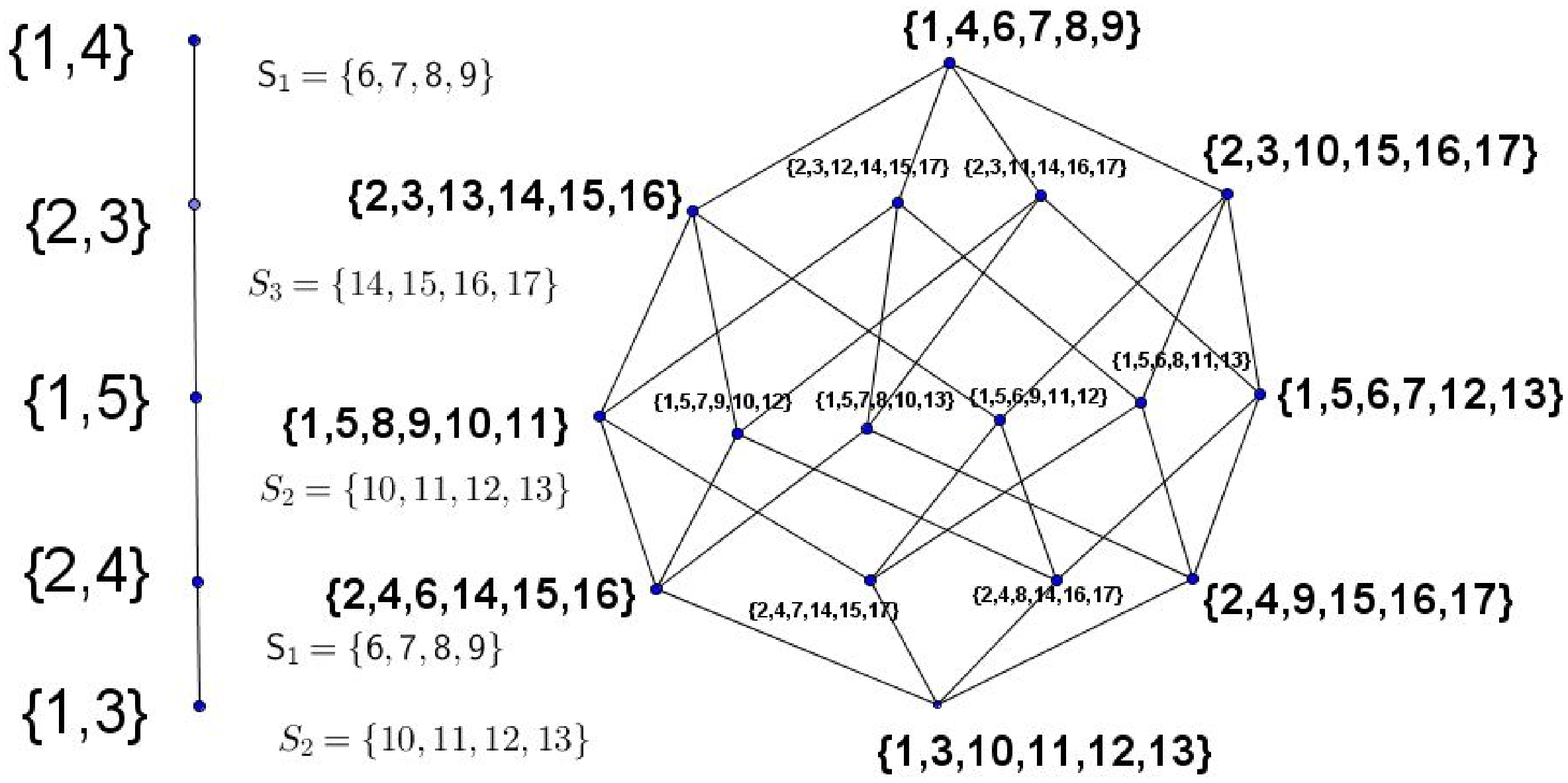}}
  \end{tabular}
  \caption{Labelling after step $1$ and final uniform labelling for $Q_{4}$}
\end{figure}
  
\section{Results on {\sc iln} }
\begin{theorem}
{\sc iln}$(P_n) = O(logn)$.
\end{theorem}
Procedure addedge may increase the cardinality of individual labels by at most one. Hence New {\sc iln}=Old {\sc iln}+ 1 (After $i^{th}$ iteration,    where $i\geq 1$) 
If we obtain $P_n$ using addedge procedure, 
then {\text {\sc iln}}$(P_{n})={\text{\sc usn}}(P_{n})/3$\\
So {\sc iln}($P_{n}$)= $O(\log{n})$ \\
The same idea is also applicable on cycles as well as wheel graph. 
\begin{theorem}
{\sc iln}($C_{n}$) = $O(logn)$.
\end{theorem}
\begin{theorem}
{\sc iln}($W_{n}$) = $O(logn)$.
\end{theorem}
\begin{theorem}
{\sc iln}$(Q_n) < n+O(logn)$
\end{theorem}
After applying Hypercube labelling algorithm, cardinality of each label will be at most n+O(logn). 
Hence individual label number is at most $n+O(logn)$
\section*{Conclusion and Future work}
We have obtained optimal value of {\sc usn} for the complement of complete graphs, complete graphs, complete bipartite graphs and upper bound for paths, cycles, matching, hypercube, wheel graph etc. Results on path, cycle, matching and wheel graphs are asymptotically tight for all three problem variants. In the future we plan to derive optimal and/or lower bound results for hypercube, harary graph etc and solve other variants of this problem like uniform labeling with minimum number of labels.

\section*{References}
\begin{enumerate}
	\item Erdos, Paul; Goodman, A. W.; Pósa, Louis (1966), "The representation of a graph by set intersections", Canadian Journal of Mathematics 18 (1): 106–112.
	\item Szpilrajn-Marczewski, E. (1945), "Sur deux propriétés des classes d'ensembles", Fund. Math. 33: 303–307.
	\item  Balakrishnan, V. K. (1997), Schaum's outline of theory and problems of graph theory, McGraw-Hill Professional, p. 40, ISBN 978-0-07-005489-9.
	\item Michael, T. S., and Thomas Quint. "Sphericity, cubicity, and edge clique covers of graphs." Discrete Applied Mathematics 154.8 (2006): 1309-1313.
	\item Lovász, L. (1968), "On covering of graphs", in Erdős, P.; Katona, G., Proceedings of the Colloquium held at Tihany, Hungary, 1966, Academic Press, pp. 231–236.
	\item Alon, Noga (1986), "Covering graphs by the minimum number of equivalence relations", Combinatorica 6 (3): 201–206.
	\item  Garey, Michael R.; Johnson, David S. (1979), Computers and Intractability: A Guide to the Theory of NP-Completeness, W. H. Freeman, ISBN 0-7167-1045-5, Problem GT59.
	\item Opsut, R. J.; Roberts, F. S. (1981), "On the fleet maintenance, mobile radio frequency, task assignment, and traffic phasing problems", in Chartrand, G.; Alavi, Y.; Goldsmith, D. L.; Lesniak-Foster, L.; Lick, D. R., The Theory and Applications of Graphs, New York: Wiley, pp. 479–492.
	\item  Scheinerman, Edward R.; Trenk, Ann N. (1999), "On the fractional intersection number of a graph", Graphs and Combinatorics 15 (3): 341–351.
	\item Bollobás, Béla, and Andrew Thomason. "Set colourings of graphs." Discrete Mathematics 25.1 (1979): 21-26.
	\item Hegde, S. M. "Set colorings of graphs." European Journal of Combinatorics 30.4 (2009): 986-995.
	\item Balister, Paul N., E. Győri, and Richard H. Schelp. "Coloring vertices and edges of a graph by nonempty subsets of a set." European Journal of Combinatorics 32.4 (2011): 533-537.
\end{enumerate}
\end{document}